\documentclass{amsart}

\usepackage{amscd,amssymb,amsmath}
\usepackage{graphics}
\usepackage[all]{xy}

\newtheorem{theorem}{Theorem}[section]

\newtheorem{lemma}[theorem]{Lemma}
\newtheorem{proposition}[theorem]{Proposition}

\theoremstyle{definition}

\newtheorem{remark}[theorem]{Remark}
\newtheorem{question}[theorem]{Question}
\newtheorem*{acknowledgement}{Acknowledgement}

\theoremstyle{remark}

\renewcommand{\labelenumi}{(\roman{enumi})}

\DeclareFontFamily{U}{wncy}{}
\DeclareFontShape{U}{wncy}{m}{n}{<->wncyr10}{}
\DeclareSymbolFont{mcy}{U}{wncy}{m}{n}
\DeclareMathSymbol{\Sh}{\mathord}{mcy}{"58}

\newcommand\mynote[1]{\marginpar{\ \\ \small \tt #1}}
\newcommand\bel[1]{{\mynote{#1}}\begin{equation}\label{#1}}
\newcommand\mylabel[1]{\label{#1}}

\newcommand{\ZZ}{\mathbb{Z}}
\newcommand{\QQ}{\mathbb{Q}}

\newcommand{\CC}{\mathbb{C}}

\newcommand{\HH}{\mathbb{H}}
\newcommand{\PP}{\mathbb{P}}

\newcommand  {\shC}     {\mathcal{C}}
\newcommand  {\shD}     {\mathcal{D}}

\newcommand  {\shF}     {\mathcal{F}}

\newcommand  {\shI}     {\mathcal{I}}

\newcommand  {\shM}     {\mathcal{M}}

\newcommand  {\shN}     {\mathcal{N}}
\newcommand  {\shL}     {\mathcal{L}}
\newcommand  {\shQ}     {\mathcal{Q}}

\newcommand  {\foM}     {\mathfrak{M}}

\newcommand  {\foS}     {\mathfrak{S}}

\newcommand  {\foU}     {\mathfrak{U}}

\newcommand  {\foX}     {\mathfrak{X}}
\newcommand  {\foY}     {\mathfrak{Y}}

\newcommand  {\ch}      {\operatorname{ch}}

\newcommand  {\ev}      {{\text{ev}}}

\newcommand  {\Grass}   {\operatorname{Grass}}

\newcommand  {\Hom}     {\operatorname{Hom}}
\newcommand  {\Hilb}    {\operatorname{Hilb}}

\renewcommand  {\ker }  {\operatorname{ker}}
\newcommand  {\Km }     {\operatorname{Km}}

\newcommand  {\length}  {\operatorname{length}}

\newcommand  {\lra}     {\longrightarrow}

\newcommand  {\Mat}     {\operatorname{Mat}}

\newcommand  {\NS}      {\operatorname{NS}}

\renewcommand{\O}       {\mathcal{O}}

\newcommand  {\Pic}     {\operatorname{Pic}}

\newcommand  {\pr}      {\operatorname{pr}}

\newcommand  {\quadand} {\quad\text{and}\quad}

\newcommand  {\ra}      {\rightarrow}

\newcommand  {\Sym}     {\operatorname{Sym}}

\def\mydate{\number\day\space\ifcase\month \or January\or February\or March\or 
April\or May\or June\or July\or
August\or September\or October\or November\or December\fi \space\number\year}

\DeclareFontFamily{U}{wncy}{}
\DeclareFontShape{U}{wncy}{m}{n}{<->wncyr10}{}
\DeclareSymbolFont{mcy}{U}{wncy}{m}{n}
\DeclareMathSymbol{\Sh}{\mathord}{mcy}{"58}

\begin{document}

\title[Periods of Enriques Manifolds]
      {Periods of Enriques Manifolds}

\author[Keiji Oguiso]{Keiji Oguiso}
\address{Department of Mathematics, Osaka University, Toyonaka 560-0043 Osaka, Japan\\ \indent and Korea Institute for Advanced Study, Hoegiro 87, Seoul, 130-722, Korea}
\curraddr{}
\email{oguiso@math.sci.osaka-u.ac.jp}

\author[Stefan Schr\"oer]{Stefan Schr\"oer}
\address{Mathematisches Institut, Heinrich-Heine-Universit\"at,
40225 D\"usseldorf, Germany}
\curraddr{}
\email{schroeer@math.uni-duesseldorf.de}


\dedicatory{To the memory of Eckart Viehweg  \\Third revised version, February 22, 2010}

\begin{abstract}
Enriques manifolds are complex spaces whose universal coverings are
hyperk\"ahler manifolds. We introduce period domains for Enriques manifolds, establish a local Torelli theorem,
and apply period maps in various situations, involving
punctual Hilbert schemes, moduli spaces of stable sheaves,
and Mukai flops.
\end{abstract}

\maketitle
\tableofcontents
\renewcommand{\labelenumi}{(\roman{enumi})}

\section*{Introduction}
In this paper we continue our study of \emph{Enriques manifolds} started in \cite{Oguiso; Schroeer 2010}.
By definition, an Enriques manifold is a connected complex space $Y$ that is not simply connected and
whose universal covering $X$ is a hyperk\"ahler manifold. The notion of  \emph{hyperk\"ahler manifolds}
was first investigated by Beauville \cite{Beauville 1983a}, and denotes simply connected compact K\"ahler manifolds
with $H^{2,0}(X)$ generated by a symplectic form. 
Such spaces were intensively studied by Huybrechts and others (see
to \cite{Gross; Huybrechts; Joyce 2003} for an introduction). Hyperk\"ahler and Enriques manifolds
are the natural generalizations of K3 and Enriques surfaces to higher dimensions.
Enriques varieties were introduced independently in \cite{Boissiere; Nieper-Wisskirchen; Sarti 2010}.

Using punctual Hilbert schemes, moduli spaces of stable sheaves, and
generalized Kummer varieties, we constructed several examples of Enriques manifolds
\cite{Oguiso; Schroeer 2010}. The basic numerical invariant for an Enriques manifold $Y$ called the \emph{index} is
the order $d\geq 2$ of its fundamental group, which is necessarily a finite cyclic group.
Most constructions yield index $d=2$ and are related to Enriques surfaces.
However, there are also examples with index $d=3,4$ coming from bielliptic surfaces.

The goal of this paper is to study \emph{periods} for Enriques manifolds, that is, linear algebra
data coming from Hodge theory, which shed some light on deformations and moduli. Throughout, we build on
the vast theory of periods for K3 surfaces, Enriques surfaces, and hyperk\"ahler manifolds.
The first main result of this paper is a \emph{Local Torelli Theorem} for 
Enriques manifolds: Roughly speaking, the base of the Kuranishi family for an Enriques manifold
is biholomorphic to some open subset of a bounded symmetric domain.
It turns out that the bounded symmetric domains in question are of type IV
for index $d=2$. In contrast, for $d\geq 3$ we have domains of type I that are biholomorphic
to complex balls.

Our notion of marking for  Enriques manifolds depends on two simple observations:
First, the fundamental group $\pi_1(Y)$ can be canonically identified with
the group of complex roots of unity $\mu_d(\CC)$, via the trace of the representation
on $H^{2,0}(X)$. Second, complex representations of   $\mu_d(\CC)$
correspond to weight decompositions $V=\bigoplus V_i$, which are indexed by the character group
$\ZZ/d\ZZ$.
Thus our period domains will be of the form
$$
\shD_L=\left\{[\sigma]\in\PP(L_{\CC,1})\mid (\sigma,\sigma)=0 \quadand (\sigma,\overline{\sigma})>0\right\},
$$
where $L_{\CC,1}$ is the weight space for the identity character  of the complexification
of a certain lattice $L$ endowed with an orthogonal representation of $G=\mu_d(C)$,
and a marking of an Enriques manifold $Y$ is an isomorphism $\phi:H^2(X,\ZZ)\ra L$,
where $H^2(X,\ZZ)$ is the \emph{Beauville--Bogomolov lattice} endowed with the canonical
representation of $G=\pi_1(Y)$.

As an application of the Local Torelli Theorem, we shall prove that any small deformation
of the known Enriques manifolds 
$$
\Hilb^n(S)/G\quadand M_H(\nu)/G \quadand \Km^n(A)/G,
$$
which come from punctual Hilbert schemes, moduli spaces of stable sheaves, and
generalized Kummer varieties, is of the same form. Note that the situation
for hyperk\"ahler manifolds is rather different.

We also show that birationally equivalent Enriques manifolds have
identical periods. Examples of birational maps are given by
\emph{Mukai flops} of $\Hilb^n(S)$, where $S$ is a K3 surface arising as
a universal covering of an Enriques surface, and the Mukai flop are  given with respect
to certain $\PP^n=\Hilb^n(C)$, where $C\subset S$ are $(-2)$-curves.
We give a detailed study of Mukai flops defined on generalized Kummer varieties
$\Km^n(A)\subset\Hilb^{n+1}(A)$ for certain abelian surfaces $A$
admiting fibrations $\varphi:A\ra E$ onto elliptic curves.
Here the Mukai flops are defined with the help of relative Hilbert schemes
$\Hilb^{n+1}(A/F)$. 
Along the way, we obtain   new examples of nonk\"ahler manifolds with trivial
canonical class that are bimeromorphic to hyperk\"ahler manifolds.

This paper is dedicated to the memory of Eckart Viehweg.
We both learned a lot from him: about moduli and many other things.

\begin{acknowledgement}
The first author is supported by JSPS Program 22340009 and by KIAS Scholar Program.
We   thank the referee for  careful reading, corrections and   suggestions.
\end{acknowledgement}

\section{Enriques manifolds and Kuranishi family}
\mylabel{Enriques}

Recall that a hyperk\"ahler manifold $X$ is a   compact complex K\"ahler manifold that
is simply connected, with $H^{2,0}(X)=H^0(X,\Omega^2_X)$ generated by a
2-form that is everywhere nondegenerate. The dimension of such manifolds is even,
and usually written as $\dim(X)=2n$.
An \emph{Enriques manifold} is a connected complex manifold $Y$ that is not simply connected,
and whose universal covering $X$ is hyperk\"ahler. Such manifolds are necessarily projective.
The trace of the representation of $G=\pi_1(Y)$
on $H^{2,0}(X)$ gives a homomorphism $G\ra\CC^\times$, which induces a canonical bijection
$G\ra\mu_d(\CC)$ with the multiplicative group of $d$-th complex roots of unity
(see \cite{Oguiso; Schroeer 2010}, Section 2, and   \cite{Beauville 1983b}, Section 4).
Throughout, we identify the groups
$$
G=\pi_1(Y)=\mu_d(\CC).
$$
The integer $d\geq 2$ is called the \emph{index} of the Enriques manifold $Y$.

Recall that the group of characters $\mu_d(\CC)\ra\CC^\times$ is cyclic of order $d$,
and contains a canonical generator, the identiy character $\zeta\mapsto\zeta$.
Throughout, we use the identification $\Hom(\mu_d(\CC),\CC^\times)=\ZZ/d\ZZ$.
A finite-dimensional complex representation   of $G$ is nothing but a finite-dimensional
complex vector space $V$ endowed with a \emph{weight decomposition} $V=\bigoplus V_i$ indexed by the characters
$i\in\ZZ/d\ZZ$. Explicitely, the \emph{weight spaces} $V_i\subset V$ is the set of vectors where each group element $\zeta\in G$ acts
via multiplication by the complex number $\zeta^i\in\CC$.  
Note that $V_0\subset V$ is the $G$-invariant subspace, and $V_1\subset V$ is the subspace where the action of
each $\zeta$ is multiplication by itself.

Now let $Y$ be an Enriques manifold of index $d\geq 2$, and
$X\ra Y$ be its universal covering, such that $X$ is a hyperk\"ahler manifold.
The fundamental group $G=\pi_1(Y)=\mu_d(\CC)$ acts on $H^1(X,\Theta_X)$,
such that we have an weight   decomposition of cohomology vector spaces
$$
H^q(X,\Theta_X)=\bigoplus_{i\in\ZZ/d\ZZ} H^q(X,\Theta_X)_i
$$
indexed by the characters of $\pi_1(Y)=\mu_d(\CC)$. Recall that $H^q(X,\Theta_X)_0$ is nothing
but the $G$-invariant part.

\begin{proposition}
\mylabel{theta invariants}
The group of global vector fields $H^0(Y,\Theta_Y)$ vanishes, and  we have 
$H^1(Y,\Theta_Y)=H^1(X,\Theta_X)_0\simeq H^{1,1}(X)_1$.
\end{proposition}

\proof
Choose a nonzero $\sigma_X\in H^{2,0}(X)$. Then $\delta\mapsto(\delta'\mapsto \sigma_X(\delta,\delta'))$
yields  an isomorphism $\Theta_X\ra\Omega^1_X$. By definition, each group element $\zeta\in G$
acts via multiplication with the complex number $\zeta\in\CC$ on $H^{2,0}(X)$.
In turn, our isomorphism induces a bijection $H^1(X,\Theta_X)_0\ra H^{1,1}(X)_1$ between weight spaces.

The projection $h:X\ra Y$ is finite and \'etale, such that the canonical map
$\Omega_Y^1\otimes_{\O_Y}h_*(\O_X)\ra h_*(\Omega^1_X)$ is bijective.
Taking duals into $h_*(\O_X)$ and using canoncial identifications, we obtain a bijection $h_*(\Theta_X)\ra \Theta_Y\otimes_{\O_Y}h_*(\O_X)$,
and the equality $H^q(Y,\Theta_Y)=H^q(X,\Theta_X)_0$ follows. Now $H^0(X,\Theta_X)=0$ ensures
$H^0(Y,\Theta_Y)=0$.
\qed

\medskip
Let $\foY\ra B$ be a \emph{Kuranishi family} of $Y=\foY_0$, that is, a deformation
of $Y=\foY_0$ that is versal and has  the property  that  $\dim H^1(Y,\Theta_Y)$   equals the embedding dimension of $0\in B$.

\begin{proposition}
\mylabel{kuranishi family}
After shrinking $B$ if necessary, the Kuranishi family $\foY\ra B$ of an Enriques manifold $Y$ of index $d\geq 2$
is universal,  the base  is smooth, and each fiber $\foY_b$ is an Enriques manifold of index $d$.
\end{proposition}

\proof
This is  a special case of general results due to Fujiki (\cite{Fujiki 1987}, Lemma 4.14)  
and Ran (\cite{Ran 1992}, Corollary 2). We recall
the arguments, since the explicit construction will be useful later. Let $\foX'\ra D'$ be the Kuranishi family
of $X=\foX'_0$. After shrinking $D'$ if necessary, we may assume that the family is universal, has
smooth base, and all its fibers are hyperk\"ahler manifolds (see \cite{Huybrechts 1999}). 
By universality, the fundamental group $G=\pi_1(Y)$ acts on this family, such that the origin $0\in D'$ is fixed.
Since the $G$-fixed locus in $\foX'$ is closed and the projection $\foX'\ra D'$ is proper, we may
also assume that $G$ acts freely on $\foX'$. It is well-known that there is a regular system of parameters $u_1,\ldots,u_r\in\O_{D',0}^\wedge$
so that the generator of $G$ acts via $u_i\mapsto e^{2\pi \sqrt{-1} n_i/d} u_i$, for certain exponents $n_i$ (see, for example, \cite{Schroeer 2004}, Lemma 5.4).
This implies that the $G$-fixed locus $D\subset D'$   is smooth of
dimension   $\dim H^1(X,\Theta_X)_0=\dim H^1(Y,\Theta_Y)$. Consider the induced family $\foX=\foX'\times_{D'}D$.
Then $G$ acts fiberwise on  $\foX$, and we obtain a family of Enriques manifolds $\foX/G\ra D$ of index $d$.
We have a commutative diagram
\begin{equation}
\label{kodaira-spencer maps}
\begin{CD}
H^1(Y,\Theta) @>>> H^1(X,\Theta)\\
@AAA @AAA\\
\Theta_D(0) @>>> \Theta_{D'}(0)
\end{CD}
\end{equation}
where the vertical maps are the Kodaira--Spencer maps. The map on the right is bijective, and
the map on the left is the induced map on $G$-invariant subspaces. Consequently, the Kodaira--Spencer map
for $\foY\ra D$ is bijective as well.  It follows that $\foY/G\ra D$ is versal, and even universal because $H^0(Y,\Theta_Y)=0$.
\qed

\section{Period domains and local Torelli}
\mylabel{Period domains}

Our next task is to define period domains and period maps for 
Enriques manifolds, in analogy to the case of Enriques surfaces (for the letter, we refer
to  \cite{Barth; Hulek; Peters; Van de Ven 2004}, Chapter VIII, Section 19). 
To this end we need a suitable notion of \emph{marking}.  Let $Y$ be an Enriques manifold and $X$ be the universal covering,
and   $H^2(X,\ZZ)$ be the \emph{Beauville--Bogomolov lattice}, which is endowed
with the primitive and integral \emph{Beauville--Bogomolov form} (see \cite{Gross; Huybrechts; Joyce 2003}, Section 23)
and an orthogonal representation of $G=\pi_1(Y)=\mu_d(\CC)$. Note that these forms and lattices are 
also called Beauville--Bogomolov--Fujiki forms and lattices.
On the complexification $H^2(X,\CC)$, we denote
by $(\sigma,\sigma')$ the induced bilinear form, such that $(\sigma,\overline{\sigma'})$ is
the induced Hermitian form. A little care has to be taken not to confuse
bilinear and   Hermitian extensions. In the following, we find it practical to say
that a nondegenerate lattice or hermitian form has signature of type $(p,*)$ if its signature
is $(p,q)$ for some integer $q\geq 0$. Our starting point is the following observation:

\begin{lemma}
\mylabel{signatures}
The lattice $H^2(X,\ZZ)$ is nondegenerate with signature of type $(3,*)$.
The  Hermitian form on the weight space $H^2(X,\CC)_1$ is nondegenerate, and has signature of type $(2,*)$
for $d=2$, and $(1,*)$ for $d\geq 3$.
\end{lemma}

\proof
According to \cite{Beauville 1983a}, Theorem 5,   the Beauville--Bogomolov lattice
 $H^2(X,\ZZ)$ is nondegenerate and has signature $(3,*)$. Since the $G$-action is orthogonal,
the eigenspace decomposition on $H^2(X,\CC)$ is orthogonal, whence the restriction of the
Beauville--Bogomolov form to each eigenspace remains nondegenerate.

In the case $d\geq 3$, the weight space $H^2(X,\CC)_1$ contains $H^{2,0}(X)$, 
but is orthogonal, with respect to the Hermitian form,  to $H^{0,2}(X)$ and the ample class coming from $Y$, 
whence has signature of type $(1,*)$.
In case $d=2$, we have $1=-1$ in the character group $\ZZ/d\ZZ$, such that the weight space contains also $H^{0,2}(X)$.
Consequently, the signature is of type $(2,*)$.
\qed

\medskip
Now let $L$ be an abstract nondegenerate  lattice with signature of type $(3,*)$, endowed with an orthogonal
representation of the cyclic group $G=\mu_d(\CC)$. We further impose the
condition that  the Hermitian form on the weight space $L_{\CC,1}$ is nondegenerate, with signature of type $(2,*)$ in case $d=2$,
and $(1,*)$ for $d\geq 3$.

An \emph{$L$-marking} for an Enriques manifold  $Y$ of index $d\geq 2$ is an equivariant isometry $\phi:H^2(X,\ZZ)\ra L$,
where $X$ is the universal covering of $Y$  
and $H^2(X,\ZZ)$ is the Beauville--Bogomolov lattice for the hyperk\"ahler manifold $X$, 
endowed with the canonical action of $G=\pi_1(Y)=\mu_d(\CC)$. 
We now define the \emph{period domain $\shD_L$ for $L$-marked Enriques manifolds} as 
$$
\shD_L=\left\{[\sigma]\in\PP(L_{\CC,1})\mid(\sigma,\sigma)=0 \quadand (\sigma,\overline{\sigma})>0\right\},
$$
where $\overline{\sigma}$ denotes complex conjugation inside the complexification $L_\CC$.
Note that for $d=2$, the weight space $L_{\CC,1}\subset L_\CC$ is invariant under complex conjugation.
On the other hand, for $d\geq 3$, each $\sigma\in L_{\CC,1}$ satisfies
$(\sigma,\sigma)=(\zeta\sigma,\zeta\sigma)=\zeta^2(\sigma,\sigma)$
for all $\zeta\in G$, whence the weight space $L_{\CC,1}\subset L_\CC$ is totally isotropic;
now the period domain is actually given by 
$$
\shD_L=\left\{[\sigma]\in\PP(L_{\CC,1})\mid(\sigma,\overline{\sigma})>0\right\}.
$$
Clearly, our period domains inside $\PP(L_{\CC,1})$ are locally closed with respect to the classical topology,
whence inherit the structure of a complex manifold. 

It turns out that $\shD_L$ is
a \emph{bounded symmetric domain}. By   results of E.\ Cartan \cite{Cartan 1935}, each bounded symmetric domain
is the product of irreducible bounded symmetric domains, and   the irreducible bounded symmetric
domains fall into six classes. The first four  are called
\emph{Cartan classical domain}, and  in Siegel's notation (\cite{Siegel 1948}, Chapter XI, \S 48) are denoted by I, II, III, IV. Recall that 
the Cartan classical domains of type $\text{\rm I}_{m,n}$ consists of complex matrices $A\in\Mat_{m\times n}(\CC)$ so that the Hermitian
matrix   $E_m - A\overline{A}^t$ is positive definite.
The Cartan classical domains of type $\text{IV}_n$ is a connected component of the set of all nonzero $z\in\CC^{n+2}$ with $z^tHz=0$ and $\bar{z}^tHz>0$, up
to nonzero scalar factors, where $H$ is a Hermitian form of signature $(2,n)$.
One should bear in mind that the symmetric bounded domain of type $\text{IV}_{2}$ is not irreducible, rather
biholomorphic to $\HH\times\HH$.

\begin{proposition}
\mylabel{classical domains}
Set $q=\dim(L_{\CC,1})$.
For $d=2$, the period domain $\shD_L$ is the disjoint union of two copies
of bounded symmetric   domains  of Type $\text{\rm IV}_{q-1}$ of dimension $q-1$.
For $d\geq 3$, the period domains $\shD_L$ are bounded symmetric    domains of
type $\text{\rm I}_{1,q-1}$, whence   biholomorphic to the complex ball of dimension $q-1$.
\end{proposition}

\proof
By our assumptions on $L$, the weight space $L_{\CC,1}$ has signature $(2,q-2)$ in case $d=2$, so the first statement holds.
Now suppose $d\geq 3$. Now $L_{\CC,1}$ has signature $(1,q-1)$, and we may identify $\shD_L$ with the set
of 
$$
\left\{(z_0:\ldots:z_{q-1})\in\PP^{q-1}\mid z_0\overline{z}_0-\sum_{i=1}^{q-1}z_i\overline{z}_i>0 \right\},
$$
which obviously coincides with the complex ball
$$
\left\{(z_1,\ldots,z_{q-1})\in\CC^{q-1}\mid \sum_{i=1}^{q-1}z_i\overline{z}_i<1\right\}.
$$
The assertion follows.
\qed

\begin{remark}
For $d=4$, such  constructions   already appeared in Kondo's study of periods for 
nonhyperelliptic curves of genus three (\cite{Kondo 2000}, \S 2).
\end{remark}

\medskip
Let $(Y,\phi)$ be an $L$-marked Enriques manifold of index $d\geq 2$, with universal covering $X$.
Let $\sigma_X\in H^{2,0}(X)$ be a nonzero form, which is unique up to scalar factors. Considered as a class
in $H^2(X,\ZZ)$, we have 
$$
(\sigma_X,\sigma_X)=0\quadand (\sigma_X,\overline{\sigma}_X)>0\quadand \sigma_X\in H^2(X,\CC)_{1}.
$$
We thus define the \emph{period point}
of our $L$-marked Enriques manifold as the induced point $[\phi(\sigma_X)]\in\shD_L$.

Now let $f:\foY\ra B$ be a flat family  of   Enriques manifolds, say  over some simply connected
complex space $B$. It follows from  Proposition \ref{kuranishi family} that each fiber $\foY_b$
is an Enriques manifold of index $d$. Moreover, the universal covering $\foX\ra\foY$
is fiber wise the universal covering, and we obtain a flat family $\foX\ra B$ of hyperk\"ahler manifolds.

Suppose we have an $L$-marking $\phi:H^2(X,\ZZ)\ra L$, where $X=\foX_0$ is the universal covering of $Y=\foY_0$.
Since the local system $R^2\phi_*\ZZ_\foX$ is constant, our $L$-marking of $Y$ \emph{uniquely} extends to
an $L$-marking 
$\phi:R^2\phi_*\ZZ_\foX\ra L_B$ of the flat family of Enriques manifolds.
In turn, we obtain a  period map
$$
p:B\lra\shD_L,\quad
b\longmapsto [\phi(\sigma_{\foX_b})].
$$
of the marked family.
Such period maps are holomorphic, according to general results of Griffiths \cite{Griffiths 1968}.
It turns out that the \emph{Local Torelli Theorem} holds:

\begin{theorem}
\mylabel{local torelli}
Let $Y$ be an $L$-marked Enriques manifold and $\foY\ra B$ be the Kuranishi family of $Y=\foY_0$.
Then   the period map $p:B\ra\shD_L$ is a local isomorphism at $0\in B$.
\end{theorem}

\proof
Since both $B$ and $\shD$ are smooth, it suffices to check that  the differential
of the period map at $0\in B$ is injective and that $B$ and $\shD$ have the same dimension.
By the Local Torelli Theorem for hyperk\"ahler manifolds (\cite{Beauville 1983a}, Theorem 5), the differential of
the period map $b\mapsto[\phi(\sigma_{\foX_b})]$ for the Kuranishi family of $X$ is bijective. 
In light of the commutative diagram (\ref{kodaira-spencer maps}), 
the differential of the period map for the Kuranishi family
of $Y$ is injective as well.

It remains to compute  vector space dimensions.
The tangent space at $0\in B$ is isomorphic to $H^1(Y,\Theta_Y)=H^1(X,\Theta_X)_0=H^{1,1}(X)_1$.
Let us first consider the case $d\geq 3$. Then $\shD_L\subset\PP(L_{\CC,1})$ is an open subset, with
respect to the classical topology, and the tangent space at the period point is 
$$
\Hom(\CC\phi(\sigma_X), L_{\CC,1}/\CC\phi(\sigma_X)) = \Hom(\CC\sigma_X,H^2(X,\CC)_1/\CC\sigma_X).
$$
The Hodge decomposition of $H^2(X,\CC)$ is invariant with respect to automorphisms of $X$, such that
$$
H^2(X,\CC)_1=H^{0,2}(X)_1\oplus H^{1,1}(X)_1\oplus H^{2,0}(X)_1. 
$$
The first summand vanishes, because $H^{0,2}(X)=H^{0,2}(X)_{-1}$ and $1\neq -1$ in the character group
$\ZZ/d\ZZ$. It follows that $H^{1,1}(X)_1$ and $\Hom(\CC\sigma_X,H^2(X,\CC)_1/\CC\sigma_X)$ have the same
dimensions.

We finally treat  the case $d=2$. Now the period domain $\shD_L$ is an open part
of a quadratic $(\sigma,\sigma)=0$ inside $\PP(L_{\CC,1})$, so the tangent space at the period point is given by
$$
\Hom(\CC\sigma_X,V/\CC\sigma_X),
$$
where $V\subset H^2(X,\CC)_1$ is the orthogonal complement of $\sigma_X\in H^2(X,\CC)$.
Clearly, we have
$$
H^2(X,\CC)_1= H^{1,1}(X)_1\oplus\CC\sigma_X\oplus\CC\overline{\sigma}_X.
$$
Taking into account  that the Beauville--Bogomolov form has 
$(\sigma_X,\sigma_X)=0$ and $(\sigma_X,\overline{\sigma}_X)>0$, the orthogonal complement in question is
$V=H^{1,1}(X)_1\oplus \CC\sigma_X$, and the argument concludes as in the preceding paragraph.
\qed

\medskip
Let $\shM_L$ be the set of isomorphism classes of $L$-marked Enriques manifolds.
Using the Local Torelli Theorem as in \cite{Gross; Huybrechts; Joyce 2003}, Definition 25.4, we conclude
that there is a unique topology and complex structure on $\shM_L$ making all the period maps defined
on the bases of the Kuranishi family holomorphic. We thus have a \emph{coarse moduli space} $\shM_L$
of $L$-marked Enriques manifolds, and the global period map
$$
p:\shM_L\lra\shD_L
$$
is \'etale. Note, however, that $\shM_L$ is not Hausdorff, as we shall see in Section \ref{Birational}.
We note in passing that the automorphism group of a marked Enriques manifold is finite,
since the same holds for hyperk\"ahler manifolds (\cite{Huybrechts 1999}, Section 9).

\section{Applications of Local Torelli}
\mylabel{Applications}

Let $S'$ be an Enriques surface. Then $G=\pi_1(S')$ is cyclic of order two,
and the universal covering $S$ is a K3 surface. Let $n\geq 1$ be an odd number.
According to \cite{Oguiso; Schroeer 2010}, Proposition 4.1, the induced $G$-action on
$X=\Hilb^n(S)$ is free, and $Y=X/G$ is an Enriques manifold of index $d=2$.
Using period maps, we now show that any small deformation of $Y$ is of the same form.

Let $\foY\ra B$ be the Kuranishi family of the Enriques manifold $Y=\foY_0$,
and denote by $\foS'\ra B'$ the Kuranishi family of the Enriques surface $S'=\foS'_0$.
We may assume that the base spaces $B$ and $B'$ are smooth and contractible.
Recall that $\dim(B')=10$. Let $\foS\ra\foS'$ be the universal covering, such that
$\foS\ra B'$ is a flat family of K3 surfaces.
The relative Hilbert scheme, or rather the relative Douady space \cite{Pourcin 1969}, 
gives a deformation $\Hilb^n(\foS/B')/G\ra B'$
of the Enriques manifold $Y$, which in turn yields a classifying map $h:B'\ra B$.

\begin{proposition}
\mylabel{versal hilbert}
The classifying map  $h:B'\ra B$ is a local isomorphism at the origin $0\in B'$.
In other words, any small deformation of the Enriques manifold $Y=\Hilb^n(S)/G$ is again of this form.
\end{proposition}

\proof
Let $L=H^2(X,\ZZ)$ be the Beauville--Bogomolov lattice, endowed with the canonical $G$-action.
Recall that Beauville \cite{Beauville 1983a} defined an  injection of Hodge structure
$$
i:H^2(S,\ZZ)\lra H^2(X,\ZZ)
$$
compatible with the $G$-action, where  the Beauville--Bogomolov form restricts to the cup product.
The cokernel is generated by  the class of the exceptional divisor of
the Hilbert--Chow map $\Hilb^n(S)\ra\Sym^n(S)$, which is $G$-invariant. It follows that we obtain an identification
$H^2(S,\CC)_1=H^2(X,\CC)_1$ of weight spaces.

Now set $L=H^2(X,\ZZ)$, and define $L'\subset L$ as the image of $i$.
In this way we obtain a marking of the Enriques manifold $Y$ and the Enriques surface $S'$.
These marking extends uniquely to markings of the families $\foY\ra\foU$ and $\foS'\ra\foU'$.
Now recall that the period domain for $L$-marked Enriques manifolds is
$$
\shD=\left\{[\sigma]\in\PP(L_{\CC,1})\mid (\sigma,\sigma)=0 \quadand (\sigma,\overline{\sigma})>0\right\}.
$$
This coincides with  the period domain for $L'$-marked Enriques surfaces  
as described in \cite{Barth; Hulek; Peters; Van de Ven 2004}, Chapter VIII, Section 19, because we have an 
equality of weight spaces $L'_{\CC,1}=L_{\CC,1}$.
Now consider the diagram
$$
\begin{xy}
\xymatrix{
B'\ar[rr]^h\ar[rd]_{p'} &       & B\ar[ld]^{p}\\
                        & \shD
}
\end{xy}
$$
where $p',p$ are the period maps for the marked families of Enriques surfaces and Enriques manifolds, respectively.
By the Local Torelli Theorems, both period maps are local isomorphisms at the origins.
It thus remains to check that the diagram is commutative. But this follows from the very definition of
the period map and the fact   that the map $i$ of local systems  sends the lines $H^{2,0}(\foS_b)$ to $H^{2,0}(\foX_b)$.
\qed

\begin{remark}
According to \cite{Fantechi 1995}, every small deformation of a punctual Hilbert scheme of an Enriques surface is of the same form, 
and our Enriques manifolds $Y=\Hilb^n(S)/G$ show the same behavior.
The situation for the hyperk\"ahler manifold $X=\Hilb^n(S)$ is different:
Its Kuranishi family has a $21$-dimensional base, whereas the Kuranishi family for the
K3 surface has only dimension $20$.
As explained in \cite{Beauville 1983a}, Theorem 6, a very general small deformation of $X$ is not isomorphic to a punctual Hilbert scheme.
\end{remark}

Keeping the previous assumptions, we now additionally assume that our Enriques surface $S'$ is general
in the sense that the $K3$ surface $S$ has the  minimal possible Picard number $\rho(S)=10$.
Let $\nu=(r,l,\chi-r)\in H^{\ev}(S,\ZZ)$ be a primitive \emph{Mukai vector}, 
with $l\in\Pic(S)$ and $\nu^2\geq 0$ and $\chi$ odd,
and $H\in \NS(S)$ a very general polarization. Then the moduli space of $X=M_H(v)$ of $H$-stable sheaves
on $S$ with Mukai vector $\nu(\shF)=\nu$ is a hyperk\"ahler manifold of dimension $\nu^2+2$.
According to \cite{Oguiso; Schroeer 2010}, Theorem 5.3, the canonical $G$-action leaves this moduli space invariant
and acts freely on it, such that $Y=X/G$ is an Enriques manifold of index $d=2$.
We shall see that any small deformation of $Y$ is of the same form.

Let $\foY\ra B$ be the Kuranishi family of the Enriques manifold $Y=\foY_0$,
and $\foS'\ra B'$ be the Kuranishi family of the Enriques surface $S'=\foS_0$.
Then we have a relative moduli space $\foM_{H}(\nu)$ for the induced family $\foS\ra B'$ of K3 surfaces,
where   $H$ now denotes a very general relative polarization.
The fiber wise $G$-action on $\foM_H(\nu)$ is free, according to \cite{Oguiso; Schroeer 2010}, Proposition 5.2,
and we obtain a flat family $\foM_H(\nu)/G\ra B'$ of Enriques manifolds.
Let $h:B'\ra B$ be the classifying map.

\begin{proposition}
\mylabel{versal moduli}
The classifying map $h:B'\ra B$ is a local isomorphism at the origin $0\in B'$.
In other words, any small deformation of the Enriques manifold $Y=M_H(\nu)/G$ is
again of this form.
\end{proposition}

\proof
Mukai (see \cite{Mukai 1988} and \cite{Mukai 1987}) defined a homomorphism
$$
\theta:\nu^\perp\lra H^2(M_H(\nu),\CC)
$$
where $\nu^\perp\subset H^\ev(S,\CC)$ denotes the orthogonal complement
with respect to the \emph{Mukai pairing}, and O'Grady \cite{O'Grady 1997} showed in full generality that
it is bijective, orthogonal, and respects the integral structure as well as the Hodge structure.
The function $\theta(\alpha)$ can be defined   on the full Mukai lattice as a K\"unneth component of 
$$
\frac{1}{\sigma}\pr_{2*}(\ch(\shQ)(1+\pr_1^*[S])\pr_1^*(\alpha))
$$
where $\pr_1:S\times M_H(\nu)\ra S$ and $\pr_2:S\times M_H(\nu)\ra M_H(\nu)$ are the
projections and $\shQ$ is a \emph{quasitautological bundle}, that is, 
a coherent sheaf on $S\times M_H(\nu)$ whose restrictions to $S\times\left\{[\shF]\right\}$ are isomorphic to $\shF^{\oplus\sigma}$ 
for some $\sigma\geq 1$, and satisfying the obvious universality property.
On the orthogonal complement $\nu^\perp$, the expression $\theta(\alpha)$ does not depend
on the choice of the quasitautological bundle. From this one infers that 
$\theta:\nu^\perp\ra H^2(M_H(\nu),\CC)$ is natural, in particular, equivariant with
respect to the canonical action of $G$.
Since the Mukai vector $\nu\in H^\ev(S,\ZZ)$ is $G$-fixed, we have $H^2(S,\CC)_1\subset\nu^\perp$,
and this yields an identification  $H^2(S,\CC)_1= H^2(M_H(\nu),\CC)_1$.
Now the argument concludes as in the previous proof.
\qed

\begin{remark}
Let $\foS\ra B$ be the flat family of K3 surfaces
induced from the Kuranishi family $\foS'\ra B$ of the Enriques surface $S'$.
According to \cite{Oguiso 2003},  every neighborhood of the origin $0\in B$ contains
points $b$ so that the fiber $\foS_b$ has Picard number $\rho>10$.
Hence there are Enriques manifolds of the form $M_H(\nu)$ arising from 
Enriques surfaces $S'$ that are more special than the general ones considered
in \cite{Oguiso; Schroeer 2010}, Section 5.
\end{remark}

\medskip
Now suppose that $S$ is a \emph{bielliptic surface}, such that  
$\omega_S\in\Pic(S)$ has order $d\in\left\{2,3,4,6\right\}$
and that the corresponding canonical covering $A$ is an abelian surface. Then $A\ra S$ is
an \'etale Galois covering, with Galois group $G=\mu_d(\CC)$. Let $n\geq 2$ be an integer with $d\mid n+1$,
and consider the generalized Kummer variety $\Km^n(A)\subset\Hilb^{n+1}(A)$ comprising those
zero cycles mapping to the origin $0\in A$ under the summation map. 
According to the results of \cite{Oguiso; Schroeer 2010}, Section 6,
with a suitable choice of the origin $0\in A$ and  for $d\neq 6$ and with one exception for $d=3$, 
the generalized Kummer variety $\Km^n(A)\subset\Hilb^{n+1}(A)$
is invariant under the canonical $G$-action  on $\Hilb^{n+1}(A)$ and the induced $G$-action on the
hyperk\"ahler manifold $X=\Km^n(A)$ is free, such that $Y=X/G$ is an Enriques manifold.
Using similiar arguments as for Proposition \ref{versal hilbert}, one shows:

\begin{proposition}
\mylabel{versal kummer}
Any small deformation of the Enriques manifold $Y=\Km^n(A)/G$ is of the same form.
\end{proposition}

\begin{remark}
One may show  
that the period domain of marked $Y=\Km^n(A)/G$  is a bounded
symmetric domain of type $I_{1,2}$ for $d=2$, whence biholomorphic to $\HH\times\HH$,
and of type $I_{1,1}$ for $d=3,4$, whence biholomorphic to $\HH$.
In both cases, it coincides with the period domain of the originial biellliptic surface $S$,
and one may use periods of elliptic curves to describe the period map explicitely.
\end{remark}

The   results of this section trigger several questions:

\begin{question}
Are the Enriques manifolds of the form $\Hilb^n(S)/G$ and $M_H(\nu)/G$ with index $d=2$ and  same dimension $2n=\nu^2+2$
deformation equivalent? This is actually true for the universal covering hyperk\"ahler manifolds by the work
of Yoshioka
(\cite{Yoshioka 2001}, Theorem 8.1, under some technical assumptions; see also 
the discussion after \cite{Markman 2010}, Theorem 2.3.), but the   techniques of deforming through
  elliptic surfcaes do not seem to carry over to an equivariant setting.

More strongly, one may ask whether each Enriques manifold of the form $M_H(\nu)/G$ is birational
to an Enriques manifold of the form $\Hilb^{n+1}(S)/G$. 
See Huybrechts work \cite{Huybrechts 1997} for results on hyperk\"ahler manifolds.
\end{question}

\begin{question}
What can be said about the image of the global period map $p:\shM_L\ra\shD_L$?
This is particularly interesting for marked Enriques manifolds of the form $\Hilb^n(S)/G$ and $M_H(\nu)/G$ of index $d=2$ 
coming from Enriques surface. In contrast to K3 surfaces, the image of the period map for Enriques surfaces is not
surjective, since it misses the classes $[\sigma]\in\PP(L_{\CC,1})$ orthogonal to some
of the $l\in L_{\CC,1}$ with $l^2=-2$.

On the other hand, the global period map for marked Enriques manifolds of the form $\Km^n(A)/G$ is surjective,
as is the case for bielliptic surfaces.
\end{question}

\section{Birational Enriques manifolds}
\mylabel{Birational}

We shall   next study birational Enriques manifolds and show that they have identical periods.
Let $Y$ be an Enriques manifold of index $d\geq 2$, and $X$ be its universal covering.
Set $L=H^2(X,\ZZ)$, and let $\phi:L\ra H^2(X,\ZZ)$ be the identity map, regarded  as an $L$-marking of $Y$.

\begin{theorem}
Let $Y'$ be another Enriques manifold that is birational to $Y$, with universal covering $X'$.
Then $Y'$ also has index $d$, and there is an $L$-marking $\phi':L\ra H^2(X',\ZZ)$
so that $(Y,\phi)$ and $(Y',\phi')$ have the same period point in
the period domain $\shD_L$.
\end{theorem}

\proof
The fundamental group is a birational invariant for smooth compact complex manifolds. 
Let $\varphi : Y \dasharrow Y'$ be a birational map, $\pi : X \to Y$ and $\pi' : X' \to Y'$ be the universal covering maps of $Y$ and $Y'$ respectively. Then $\varphi \circ \pi$ is a rational map from $X$ to $Y'$.  Let $\nu : Z \to X$ be a Hironaka's resolution of indeterminacy of $\varphi \circ \pi$. Similarly, we choose a Hironaka's resolution of indeterminacy $\nu' : Z' \to X'$ of the rational map $\varphi^{-1} \circ \pi'$. Since $Z$ is smooth and birational to $X$, it follows that $Z$ is simply connected. The same is true for $Z'$. Thus, the morphism $\varphi \circ \pi \circ \nu : Z \to Y'$ can be lifted to a morphism to $X'$, say $\varrho : Z \to X'$. Similarly, 
we have a morphism $\varrho' : Z' \to X$ which is a lift of $\varphi^{-1} \circ \pi' \circ \nu'$. By (1), $\varrho$ and $\varrho'$ are both birational. Thus 
$$
f := \varrho \circ \nu^{-1} : X \dasharrow X'
$$ 
is a birational map. Note that $f$ is isomorphic in codimension $1$. This is because $K_{X}$ and $K_{X'}$ are both trivial. Thus, we can naturally pull back 
the $2$-form on $X'$ to $X$. Thus, we can choose a generator $\sigma_{X}$ 
of $H^2(X, \Omega_{X}^2)$ and a generator $\sigma_{X'}$ 
of $H^2(X', \Omega_{X'}^2)$ such that $\sigma_{X} = f^{*} \sigma_{X'}$. Moreover, by 
\cite{Gross; Huybrechts; Joyce 2003}, Proposition 25.14 (see also Pages 213--214),  
$f$ naturally induces a Hodge isometry with respect to the Beauville--Bogomolov form: 
$$
f^* : H^2(X', \ZZ) \simeq H^2(X, \ZZ)\,\, .
$$
On the other hand, by the construction, we have $\varphi \circ \pi = \pi' \circ f$. 
Thus, $f \circ g = g \circ f$. Here $g$ is a generator of $\pi_1(Y)=\mu_d(\CC)=\pi_1(Y')$.
Hence $f^*$ induces a bijection of weight spaces:
$$
H^2(X', \CC)_1 \simeq H^2(X, \CC)_1.
$$
Hence, by $f^{*}$, the period of $Y$ and the period of $Y'$ become the same point.
\qed

\medskip
To give examples of birationally equivalent Enriques manifolds, we have to recall
certain birational transformations for hyperk\"ahler manifolds.
Let $X$ be a hyperk\"ahler manifold of dimension $\dim(X)=2n$, with $n\geq 2$.
Suppose there is a closed subspace $P\subset X$ with $P\simeq\PP^n$.
As described in \cite{Huybrechts 1997}, the \emph{Mukai flop} $\check{X}$ of $X$
with respect to $P\subset X$ is defined via a commutative diagram
\begin{equation}
\label{flop diagram}
\begin{xy}
\xymatrix{
         & \tilde{X}\ar[dl]\ar[dr]\\
X\ar[dr] &                          & \check{X}\ar[dl]\\
         & \bar{X}
}
\end{xy}
\end{equation}
Here $\tilde{X}\ra X$ is the blowing-up with center $P\subset X$. Its exceptional
divisor is isomorphic to the incidence scheme $E\subset P\times\check{P}$, where
$\check{P}$ denotes the dual projective space, and $\tilde{X}\ra X$ contracts $E$
along the first projection of $P\times\check{P}$. 
The morphism $\tilde{X}\ra \check{X}$ is defined as the contraction of $E$
along the second projection, which is also the blowing-up of $\check{P}\subset\check{X}$. 
Moreover, the morphisms $X\ra\bar{X}$
and $\check{X}\ra\bar{X}$ are the contractions of $P$ and $\check{P}$, respectively.
Note that the Mukai flop is simply connected and $H^{2,0}(\check{X})$ is generated
by a symplectic form, such that $\check{X}$ is hyperk\"ahler if and only if
it is K\"ahler.
Given several disjoint copies $P_1,\ldots,P_d\subset X$, we may also perform a Mukai flop $\check{X}$
with respect to  $P_1\cup\ldots\cup P_d\subset X$ simultaneously.

Now let $X$ be the universal covering of an Enriques manifold $Y$ of index $d\geq 2$,
and suppose there is a copy $Q\subset Y$  of $\PP^n$.
Since $X\ra Y$ is a local isomorphism on a small neighborhood of $P_i\subset X$ with respect to
the classical topology, we obtain a Mukai flop $\check{Y}$ of $Y$ with respect to $Q\subset Y$,
whose universal covering is the Mukai flop $\check{X}$ of $X$ with respect to $P_1\cup\ldots\cup P_d\subset X$.
Note, however, that it is in general not so easy  to determine whether or not the Mukai  flops are K\"ahler,
such that $\check{Y}$ is indeed an Enriques manifold.

Here is an example: Let $S'$ be an Enriques surface, with universal covering $S$,
and  $n\geq 0$ be an odd number. Then the induced action of $G=\pi_1(S')$ on the
punctual Hilbert scheme $X=\Hilb^n(S)$ is free, such that $Y=X/G$ is an Enriques manifold
of index $d=2$.
Suppose furthermore that $S'$ is nodal, that is, there is a curve $C'\subset S'$ of arithmetic genus zero.
Then $C'\simeq \PP^1$ and $C'^2=-2$, that is, $C'\subset S'$ is a $(-2)$-curve.
The preimage $C_1\cup C_2\subset S$ is a   union of two disjoint $(-2)$-curves.
In turn, we obtain two disjoint copies
$$
P_i=\Hilb^n(C_i)\subset \Hilb^n(S)=X,\quad i=1,2
$$
of $\PP^n=\Sym^n(\PP^1)=\Hilb^n(\PP^1)$ inside the hyperk\"ahler manifold, which are interchanged by the $G$-action.
Set  $Q=(P_1\cup P_2)/G\subset Y$.

\begin{proposition}
Assumptions as above. Then
the Mukai flop $\check{Y}$ of $Y$ with respect to $Q\subset Y$ is projective, whence an Enriques manifold
birational to $Y$.
\end{proposition}

\proof
This is a variation of  an argument of Debarre \cite{Debarre 1984}, where maps
to Grassmannians are exploited.
To carry  it out, we have to  verify that there is some $\shL'\in \Pic(S')$ with the following properties:
$\shL'\cdot C'=1$ and both $\shL'$ and its preimage $\shL\in\Pic(S)$ are very ample.
In other words, $C'$ and $C_1,C_2$ become \emph{lines} under suitable embeddings into  
projective spaces. To see this, consider the 
contraction $S'\ra \bar{S}'$  of the $(-2)$-curve $C'$. Then the proper normal surface $\bar{S}'$ is projective.
Let $D_1$ be the pullback of some ample line bundle.
Since the intersection form on $\NS(S')$ is unimodular, there is a  divisor $D_2$ with 
$D_2\cdot C'=1$. Consider the invertible sheaf $\shL'=\O_{S'}(nD_1+D_2)$. Then
$\shL'(-C'-K_{S'})$ is  relatively ample over $\bar{S}'$, whence ample for $n\gg 0$, 
such that $H^1(S',\shL'(-C'))=0$ by Kodaira Vanishing.
Arguing similarly on $S$, we easily infer that $\shL'$ has the desired
properties for $n\gg 0$. Increasing $n$ if necessary, we furthermore achieve that
$C'\subset S'$ and $C_1,C_2\subset S$ are the only lines with respect to $\shL'$ and $\shL$,
respectively. Using maps to Grassmannians as 
in \cite{Debarre 1984}, Section 3.2, we see that the Mukai flop $\check{X}$  is indeed projective,
and this then   holds for  $\check{Y}$ as well.
\qed

\begin{remark}
As discussed in \cite{Huybrechts 1999}, Section 2.4, the existence of such Mukai flops
implies that the coarse moduli space of $L$-marked Enriques manifolds $\shM_L$ is
not Hausdorff: any two birationally equivalent Enriques manifolds   give rise
to \emph{nonseparated points} of the moduli space.

Of course, it is another matter whether or not a Mukai flop
$Y\dashrightarrow\hat{Y}$ yields Enriques manifolds that are  nonisomorphic. 
Examples of bimeromorphic yet nonisomorphic hyperk\"ahler manifolds
were constructed by Debarre \cite{Debarre 1984}.  Namikawa \cite{Namikawa 2002} even found
hyperk\"ahler manifolds having the same periods that are not bimeromorphic. Both constructions, however,
use nonprojectivity in an essential way, and apparently do not apply   to Enriques manifolds.
\end{remark}

\section{Rational maps to Grassmannians}
\mylabel{Rational maps}

Recall that the geometry of the punctual Hilbert scheme $\Hilb^{n+1}(E)=\Sym^{n+1}(E)$ of an elliptic curve $E$
is very simple: The canonical map
$$
\Hilb^{n+1}(E)\lra\Pic^{n+1}(E),\quad [Z]\longmapsto\O_E(Z)
$$
is a fibration whose fiber $\Hilb^{n+1}_\shN(E)$ over an invertible sheaf $\shN\in\Pic^{n+1}(E)$ 
is isomorphic to the projectivization of $H^0(E,\shN)$.
Whence $\Hilb^{n+1}(E)\ra\Pic^{n+1}(E)$ is a $\PP^n$-bundle. Moreover, we have a commutative
diagram
$$
\begin{xy}
\xymatrix{
  & \Hilb^{n+1}(E)\ar[dl]_+\ar[dr]\\
E &                              & \Pic^{n+1}(E)\ar[ll],
}
\end{xy}
$$
where the map on the left is the composition of   Hilbert--Chow   addition map,
and the horizontal map makes this diagram commutative. The latter is bijective (and not just an isogeny), because both
diagonal arrows have connected fibers. Note that it sends the invertible sheaf associated to the divisor $(n+1)0\subset E$ to
the origin $0\in E$.

Beauville \cite{Beauville 1983b} and Debarre \cite{Debarre 1984} exploited that  maps to Grassmannians  
yield interesting contractions of punctual Hilbert schemes for K3 surfaces.
We now continue this line of thought in the following situation:
Let $A$ be an abelian surface, and consider the Hilbert scheme
of points $\Hilb^{n+1}(A)$ for some fixed integer $n\geq 1$.
Given an ample $\shL\in\Pic(A)$, we obtain for each closed
subscheme $Z\subset A$ of length $n+1$  a restriction map
$$
H^0(A,\shL)\lra H^0(Z,\shL_Z),
$$
where the vector space $H^0(Z,\shL_Z)$ is $(n+1)$-dimensional.
By the Kodaira Vanishing Theorem,  $H^0(A,\shL)$ has 
dimension $c_1^2(\shL)/2$ and $H^1(A,\shL)=0$. Whence  the restriction map is surjective if and only
if $H^1(A,\shI_Z\otimes\shL)=0$, where $\shI_Z\subset\O_A$ is the ideal of
$Z\subset A$. The set of such $[Z]\in\Hilb^{n+1}(A)$ is open an non-empty, and we get a rational map
$$
r=r_\shL:\Hilb^{n+1}(A)\dashrightarrow\Grass(V,n+1),\quad 
[Z]\longmapsto (H^0(A,\shL)\ra H^0(Z,\shL_Z))
$$
into the Grassmannian of $(n+1)$-dimensional quotients of $V=H^0(A,\shL)$. We say that
this rational map \emph{is defined} at $[Z]\in\Hilb^{n+1}(A)$ if $H^1(A,\shI_Z\otimes\shL)=0$.
Now the basic observation is:

\begin{lemma}
\mylabel{contracts}
Suppose $E\subset A$ is an elliptic curve and $\shN\in\Pic^{n+1}(E)$
satisfies $\shN\neq\shL_E$ and  $\shL\cdot E=n+1$ and $H^1(A,\shL(-E))=0$.
Then the rational map $r:\Hilb^{n+1}(A)\dashrightarrow \Grass(V,n+1)$ is defined
in a neighborhood of $\Hilb^{n+1}_\shN(E)$, and maps
$\Hilb^{n+1}_{\shN}(E)$ to a point.
\end{lemma}

\proof
The exact sequence
$$
H^0(A,\shL)\lra H^0(E,\shL_E)\lra H^1(A,\shL(-E))
$$
shows that the restriction map $H^0(A,\shL)\ra H^0(E,\shL_E)$ is surjective.
To check that the rational map is defined in a neighborhood of $\Hilb^{n+1}_\shN(E)$,
it suffices to verify that for any closed subscheme $Z\subset E$ of length $n+1$ with $\shN\simeq\O_E(Z)$,
the restriction map $H^0(E,\shL_E)\ra H^0(Z,\shL_Z)$ is surjective.
Indeed,   the outer terms in the long exact sequence
$$
H^0(E,\shL_E(-Z))\lra H^0(E,\shL_E)\lra H^0(Z,\shL_Z)\lra H^1(E,\shL_E(-Z))
$$
vanish, because the invertible sheaf $\shL_E(-Z)$ has degree zero but is nontrivial.
Thus the rational map is defined at the point $[Z]$.
Furthermore, the restriction map $H^0(A,\shL)\ra H^0(Z,\shL_Z)$ factors over the
bijection $H^0(E,\shL_E)\lra H^0(Z,\shL_Z)$, which means that the image $r([Z])$ does
not depend on the point $[Z]\in\Hilb^{n+1}_{\shN}(E)$.
\qed

\medskip
For the rest of this section, we assume that our abelian surface $A$ is endowed 
with a homomorphism 
$\varphi:A\ra F$ onto an elliptic curve $F$, such that its fibers are elliptic curves. Fix an integer $n\geq 1$,
and consider the inclusion of the relative into the absolute
Hilbert scheme   $\Hilb^{n+1}(A/F)\subset\Hilb^{n+1}(A)$.
The relative Hilbert scheme comes along with the structure map $\Hilb^{n+1}(A/F)\ra F$, which factors over the
canonical map 
$$
\Hilb^{n+1}(A/F)\lra\Pic^{n+1}(A/F),
$$
and the latter is a $\PP^n$-bundle.  Throughout, we denote by $\Pic^{n+1}(A/F)\subset\Pic(A)$ the subset of invertible sheaves $\shL$
that have degree $n+1$ on the fibers of $\varphi:A\ra F$, and by $\Hilb^{n+1}_\shL(A/F)$ the family
of zero cycles $Z$ on fibers $E=\varphi^{-1}(f)$ with $\O_E(Z)\simeq\shL_E$.

\begin{proposition}
\mylabel{not defined}
For each $\shL\in\Pic^{n+1}(A/F)$, the rational map $r:\Hilb^{n+1}(A)\dashrightarrow\Grass(V,n+1)$
is not defined on the closed subset $\Hilb^{n+1}_\shL(A/F)$.
\end{proposition}

\proof
Let $E=\varphi^{-1}(f)$ be a fiber and $Z\subset E$ be a divisor of length $n+1$ with $\O_E(Z)\simeq\shL_E$.
The restriction map $H^0(A,\shL)\lra H^0(Z,\shL_Z)$ factors over the map on the left of the exact sequence
$$
H^0(E,\shL_E)\lra H^0(Z,\shL_Z)\lra H^1(E,\shL(-Z))\lra 0.
$$
Using that  $H^1(E,\shL(-Z))\neq 0$, we conclude that the restriction map is not surjective.
\qed

\begin{proposition}
\mylabel{defined}
For  all ample  $\shL'\in\Pic^{n+1}(A/F)$ and each   $\shN\in\Pic(F)$
of degree $\deg(\shN)\geq n^2+1$, the rational map $r:\Hilb^{n+1}(A)\dashrightarrow\Grass(V,n+1)$ given by  $\shL=\shL'\otimes \varphi^*(\shN)$
is defined on the complement of $\Hilb^{n+1}_{\shL}(A/F)$.
\end{proposition}

\proof
Let $Z\subset A$ be a subscheme of length $n+1$ with $[Z]\not\in\Hilb^{n+1}_{\shL}(A/F)$.
First, suppose that $Z\subset E=\varphi^{-1}(f)$ is contained in a fiber, and $\O_E(Z)\not\simeq\shL_E$.
Since $\shL=\shL'\otimes\varphi^*(\shN(-f))$ is ample, we have $H^1(A,\shL(-E))=0$,
whence Lemma \ref{contracts} tells us that the rational map is defined at $[Z]$.

Now assume that $Z$ is not contained in any fiber of $\varphi:A\ra F$.
We have do distinguish two cases:
To start with, suppose that $Z$ is reducible, and decompose $Z=Z_1+\ldots +Z_r$, $2\leq r\leq n$   into
parts whose support is  contained in pairwise different fibers $E_i=\varphi^{-1}(f_i)$.
The sheaf $\shL(-n(E_1+\ldots+E_r))=\shL'\otimes\varphi^*(\shN(-n(f_1+\ldots+f_r)))$ is ample, whence the 
term on the right in the exact sequence
$$
H^0(A,\shL)\lra \bigoplus_{i=1}^r H^0(nE_i,\shL_{nE_i})\lra H^1(A,\shL(-n(E_1+\ldots+E_r)))
$$
vanishes. Thus it suffices to check that for each zero cycle $W\subset A$ of length $\leq n$
whose support is contained in a fiber $E=\varphi^{-1}(f)$, the restriction map
$H^0(nE,\shL)\ra H^0(W,\shL_W)$ is surjective. To this end,   set $W_i=W\cap iE$, $1\leq i\leq n$.
Then $W_1\subset W_2\subset\ldots\subset W_n$ is a sequence of zero cycles with $W_1\subset E$ and $W_n=W$.
Clearly, the restriction map $H^0(E,\shL_E)\ra H^0(W_1,\shL_{W_1})$ is surjective, because $\deg(W_1)<\deg(\shL_E)$.
Let $\shI\subset\O_{W_{i+1}}$ be the ideal of $W_i\subset W_{i+1}$. The commutative diagram
$$
\begin{CD}
0 @>>> \O_E(-E) @>>> \O_{(i+1)E} @>>> \O_{iE} @>>> 0\\
&& @VVV @VVV @VVV\\
0 @>>> \shI @>>> \O_{W_{i+1}} @>>> \O_{W_i} @>>> 0
\end{CD}
$$
shows that $\shI$ is isomorphic to $\O_D$ for some divisor $D\subset E$ of length $\leq n$.
Tensoring with $\shL$, using $\O_E(-E)\simeq\O_E$ and taking cohomology, we obtain a commutative diagram with exact rows
$$
\begin{CD}
0 @>>> H^0(E,\shL_E)@>>> H^0((i+1)E,\shL_{(i+1)E}) @>>> H^0(iE,\shL_{iE}) @>>> 0\\
&& @VVV @VVV @VVV\\
0 @>>> H^0(D,\shL_D) @>>> H^0(W_{i+1},\shL_{W_{i+1}}) @>>> H^0(W_i,\shL_{W_i}) @>>> 0.
\end{CD}
$$
The vertical map on the left is surjective because $\deg(D)<\deg(\shL_E)$, and the vertical map
on the right is surjective by induction. We conclude that  $H^0(nE,\shL)\ra H^0(W,\shL_W)$ is surjective.

It remains to treat the case that $Z\subset A$ is irreducible, say with support contained in $E=\varphi^{-1}(f)$,
but not contained in $E$ as a subscheme.  We then argue as in the preceding paragraph, taking into account
that $Z_1=Z\cap E$ has length $\leq n$.
\qed

\begin{proposition}
\mylabel{injective}
For  all ample  $\shL'\in\Pic^{n+1}(A/F)$ and each   $\shN\in\Pic(F)$
of degree $\deg(\shN)\geq n^2+2$, the rational map $r:\Hilb^{n+1}(A)\dashrightarrow\Grass(V,n+1)$ given by  $\shL=\shL'\otimes \varphi^*(\shN)$
is injective on the complement of $\Hilb^{n+1}(A/F)$.
\end{proposition}

\proof
Let $Z,Z'\subset A$ be two different zero cycles of length $n+1$, none of them contained in fibers of $\varphi:A\ra F$.
We have to find some section $s\in H^0(A,\shL)$ vanishing on one but not on  both   cycles,
for then the kernels of the two surjections $H^0(A,\shL)\ra H^0(Z,\shL_Z)$ and $H^0(A,\shL)\ra H^0(Z',\shL_{Z'})$
are different.

Suppose that we can find a zero cycle $Z\subset W\subset Z\cup Z'$ with $\length(W)=n+2$
so that the intersection $W\cap\varphi^{-1}(f)$ with every fiber has length $\leq n$.
Arguing as in the preceding proof, one can infer that the restriction map
$H^0(A,\shL)\ra H^0(W,\shL_W)$ is surjective, and we are done.

It remains to treat the case that there is no such zero cycle in neither $Z\subset Z\cup Z'$
nor  $Z'\subset Z\cup Z'$. Then it easily follows that some fiber $E=\varphi^{-1}(f)$
intersects $Z$ in length $n$, and this fiber intersects $Z'$ in length $n$ as well,
and $Z,Z'$ are supported on $E$.
We are done if $W_1=(Z\cup Z')\cap E$ has length $\geq n+2$, because a nonzero section of $\shL_E$
cannot vanish on a subscheme of length greater than $\deg(\shL_E)=n+1$.
It remains to treat the case that $W_1$ has length $n+1$. 
We then argue on the infinitesimal neighborhood $2E$ as in the preceding proof.
Details are left to the reader.
\qed

\section{Mukai flops of generalized Kummer varieties}
\mylabel{Mukai flops}

Let $X$ be a hyperk\"ahler manifold of dimension  $\dim(X)=2n$, this time with $n\geq 2$.
Recall that for each closed  subspace $P\subset X$ with $P\simeq\PP^n$,
the diagram 
\begin{equation}
\begin{xy}
\xymatrix{
         & \tilde{X}\ar[dl]\ar[dr]\\
X\ar[dr] &                          & \check{X}\ar[dl]\\
         & \bar{X}
}
\end{xy}
\end{equation}
defines the Mukai flop $\check{X}$ of $X$ with respect to $P\subset X$,
and that the Mukai flop $\check{X}$ is hyperk\"ahler provided that it is K\"ahler.
The following two results, important in the sequel, are well-known:

\begin{lemma}
\mylabel{projective flop}
Assume that the Mukai flop $\check{X}$ is K\"ahler. 
Then $\check{X}$ is projective if and only if $X$ is projective.
\end{lemma}

\proof
If $X$ is projective, then $\check{X}$ is Moishezon. Being K\"ahler,
it must be projective. The reverse implication holds by symmetry.
\qed

\medskip
Now let $P,P'\subset X$ be two copies of $\PP^n$. We say that
$P,P'\subset X$ are \emph{numerically equivalent} if $\deg(\shL_P)=\deg(\shL_{P'})$
for all $\shL\in\Pic(X)$.

\begin{proposition}
\mylabel{nonprojective flop}
Assume that $P,P'\subset X$ are two disjoint copies of $\PP^n$ that are numerically
equivalent. Then the Mukai flop $\check{X}$ of $X$ with respect to $P\subset X$
is not projective. If $X$ is projective, $\check{X}$ is not even K\"ahler.
\end{proposition}

\proof
Seeking a contradiction, we assume that $\check{X}$ is projective.
Then $X$ is projective as well, by Proposition \ref{projective flop}.
Choose an ample $\check{\shL}\in\Pic(\check{X})$, and let $\shL\in\Pic(X)$ be its strict transform. Then
$$
\deg(\shL_P)=\deg(\shL_{P'})=\deg(\check{\shL}_{P'})>0,
$$
where we regard $P'\subset X\smallsetminus P=\check{X}\smallsetminus\check{P}$ also as a closed subspace
of the Mukai flop. It follows that $\shL$ is relatively ample for the contraction $h:X\ra\bar{X}$, and
obviously $\check{\shL}$ is relatively ample for $\check{h}:\check{X}\ra\bar{X}$.
Using that $P\subset X$ and $\check{P}\subset\check{X}$ have codimension $n\geq 2$, the
equality $X\smallsetminus P=\check{X}\smallsetminus\check{P}$ induces an identification
$$
\bigoplus_{i\geq 0} h_*(\shL^{\otimes i}) = \bigoplus_{i\geq 0} \check{h}_*(\check{\shL}^{\otimes i})
$$
of $\O_{\bar{X}}$-algebras. Taking the relative homogeneous spectra, 
we see that the rational map $X\dashrightarrow\check{X}$
extends to a $\bar{X}$-isomorphism $X\ra \check{X}$. In turn, the pullbacks of $\shL$ and $\check{\shL}$
to $\tilde{X}$ become isomorphic, which   gives a contradiction.

Finally, suppose that $X$ is projective. If $\check{X}$ is K\"ahler, then it is also projective.
But Proposition   \ref{nonprojective flop} tells us that it is not projective.
\qed

%
%

\medskip
Now let $A$ be again an 2-dimensional complex torus endowed with a homomorphism $\varphi:A\ra F$ onto some
elliptic curve. The punctual Hilbert scheme $\Hilb^{n+1}(A)$ contains the generalized Kummer variety $\Km^n(A)$,
which is defined as the cartesian diagram
$$
\begin{CD}
\Km^n(A) @>>> \Hilb^{n+1}(A)\\
@VVV         @VV+V\\
0 @>>> A,
\end{CD}
$$
and the relative Hilbert scheme $\Hilb^{n+1}(A/F)$.

\begin{proposition}
\mylabel{kummer intersection}
The intersection $\Km^n(A)\cap \Hilb^{n+1}(A/F)$ inside
$\Hilb^{n+1}(A)$ is the disjoint union of $(n+1)^2$ copies
of $\PP^n$.
\end{proposition}

\proof
Let $Z\subset E=\varphi^{-1}(f)$ be a zero cycle of length $n+1$ contained
in a fiber, and write $Z=\sum_{i=1}^{n+1}z_i$ with $z_i\in E$.
Applying $\varphi:A\ra F$ to the sum of $Z$ in $A$, we obtain $(n+1)f\in F$.
Thus $[Z]\in\Km^n(A/F)$ implies that $f\in F[n+1]$, hence there are only $(n+1)^2$ possibilities for $f$.

Now suppose that $f\in F[n+1]$, and set $E=\varphi^{-1}(f)$. Consider the map
$$
h_f:\Hilb^{n+1}(E)\lra \ker(\varphi)
$$
that sends a zero cycle $Z\subset E$ into its sum in $A$. Clearly,
$\Km^n(A)\cap\Hilb^{n+1}(A/F)$ is the disjoint union of the preimages $h_f^{-1}(0)$, whence
it is the disjoint union of $(n+1)^2$ copies of $\PP^n$.
\qed

\begin{proposition}
\mylabel{numerically equivalent}
The $(n+1)^2$ components of  $\Km^n(A)\cap\Hilb^{n+1}(A/F)\subset\Km^n(A)$ are pairwise numerically equivalent.
\end{proposition}

\proof
According to \cite{Beauville 1983a}, Section 7, the restriction map
$$
H^2(\Hilb^{n+1}(A),\QQ)\ra H^2(\Km^n(A),\QQ)
$$ 
is surjective. It follows that each invertible sheaf on $\Km^n(A)$ has a multiple that
is the restriction of an invertible sheaf on $\Hilb^{n+1}(A)$. The latter
have constant degree on the fibers of the flat family $\Hilb^{n+1}(A/F)\ra\Pic^{n+1}(A/F)$,
and our  $P_i\subset\Km^n(A)$ are fibers of this family.
\qed

\medskip
Let us write $P_i\subset\Km^n(A)$, $1\leq i\leq (n+1)^2$
for these copies of $\PP^n$, and  $P_I=\cup_{i\in I} P_i$ for the   disjoint union
for a given subset  $I\subset\left\{1,\ldots,(n+1)^2\right\}$.

\begin{theorem}
\mylabel{multiple nonkaehler}
Let  $I\subsetneqq\left\{1,\ldots,(n+1)^2\right\}$ be a proper subset. Then the Mukai flop $\hat{X}_I$
of $X=\Km^n(A)$ along $P_I\subset X$ is not K\"ahler.
\end{theorem}

\proof
This follows from Proposition \ref{numerically equivalent} and Proposition \ref{nonprojective flop}.
\qed

\begin{remark}
Yoshioka  constructed first examples of nonk\"ahler manifolds with trivial canonical class
bimeromorphic to hyperk\"ahler manifolds, 
using Mukai flops of moduli spaces of stable sheaves on abelian surfaces (\cite{Yoshioka 2001}, Section 4.4).
Note also that Guan \cite{Guan 1995} has used primary Kodaira surface, which are not in class $\shC$, 
to construct compact symplectic manifolds that are not in class $\shC$.
\end{remark}

On the other hand, the results of the previous section give a projectivity statement:

\begin{theorem}
\mylabel{multiple projective}
Let $\hat{X}$ be the simultaneous Mukai flop
of $X$ with respect to  the full intersection $\Km^n(A)\cap\Hilb^{n+1}(A/F)=P_1\cup\ldots\cup P_{(n+1)^2}$. 
If the homomorphism $\varphi:A\ra F$ admits a section,
then the Mukai flop $\hat{X}$ is projective.
\end{theorem}

\proof
Write $A=E\times F$ with $E=\ker(\varphi)$. Choose some divisor $D_1\subset E$
of degree $n+1$ not linearly equivalent to $(n+1)0$, and some $D_2\subset E$ of degree $\geq n^2+2$.
Consider the invertible sheaf $\shL=\O_A(D_1\times F + E\times D_2)$ and the rational map
$$
r:\Hilb^{n+1}(A)\dashrightarrow \Grass(V,n+1),\quad Z\longmapsto (H^0(A,\shL)\ra H^0(Z,\shL_Z))
$$
studied in the previous section. According to Proposition \ref{defined}, this rational map is defined
on some neighborhood of $\Km^n(A)\subset\Hilb^{n+1}(A)$.
Furthermore, it sends $P_1,\ldots,P_{(n+1)^2}$ to points, and is injective on the complement,
by Lemma \ref{contracts} and Proposition \ref{injective}. Hence the Stein factorization 
of $r:\Km^n(A)\ra\Grass(V,n)$ factors over the contraction $X\ra\bar{X}$, such that
$\bar{X}$ is projective. Using that $\tilde{X}$ is projective, one easily infers
that the resolution $\check{X}\ra\bar{X}$ is relatively projective. The upshot is
that $\check{X}$ is projective.
\qed

\medskip
Now let $A=E\times F$ be a product of two elliptic curves, on which $G=\mu_2(\CC)$
acts freely via the involution
$$
(e,f)\longmapsto (e+1/2,-f+z)
$$
as explained in \cite{Oguiso; Schroeer 2010}. Suppose $n\geq 2$ is an odd integer,
and $z\in F$ is a point with $(n+1)z=0$ but $(n+1)/2\cdot z\neq 0$. According to loc.\ cit.,
Theorem 6.4, the $G$-action on $\Hilb^{n+1}(A)$ leaves the subset $\Km^n(A) $ invariant and acts
freely there, such that $Y=\Km^n(A)/G$ is an Enriques manifold of index $d=2$.
The group $G$ permutes  the $(n+1)^2$ copies $P_i$ of $\PP^n$ inside $X=\Hilb^{n+1}(A)$ considered above,
whence the Mukai flop $\bar{X}$ with respect to any $G$-invariant union $P_I$ induces a Mukai flop
$\bar{Y}$ of our Enriques manifold $Y$.
The upshot is:

\begin{theorem}
There are Enriques manifolds $Y$ admitting   Mukai flops that are nonk\"ahler
(whence not Enriques manifolds in the strict sense, rather  ``nonk\"ahler Enriques manifolds''),
and other Mukai flops $\hat{Y}$ that are K\"ahler (so again true Enriques manifolds).
\end{theorem}

\begin{question}
Is such $\hat{Y}$ isomorphic to $Y$ as an abstract complex space?
\end{question}



\end{document}